\documentclass[12pt]{article}
\usepackage[pctex32]{graphics}
\usepackage{amsthm,amsmath,amssymb,amscd,verbatim,epsfig}
\usepackage{amssymb}
\usepackage{graphicx}
\newcommand{\beq}{\begin{equation}}
\newcommand{\eeq}{\end{equation}}

\date{}

\newcommand{\al}{\alpha}

\newcommand{\Ld}{\Lambda}
\newcommand{\ga}{\gamma}
\newcommand{\be}{\beta}

\newcommand{\si}{\sigma}
\newcommand{\Si}{\Sigma}

\newcommand{\dl}{\delta}

\newcommand{\iy}{\infty}

\newcommand{\ra}{\rightarrow}

\newcommand{\supp}{\supset}

\newcommand{\sq}{$\square$}

\begin{document}
\title{On\ the\ nature\ of\ chaos}
\author{Bau-Sen Du \\ [.5cm]
Institute of Mathematics \\
Academia Sinica \\
Taipei 11529, Taiwan \\
dubs@math.sinica.edu.tw \\}
\maketitle
\begin{abstract}
Based on newly discovered properties of the shift map (Theorem 1), we believe that chaos should involve not only nearby points can diverge apart but also faraway points can get close to each other.  Therefore, we propose to call a continuous map $f$ from an infinite compact metric space $(X, d)$ to itself chaotic if there exists a positive number $\delta$ such that for any point $x$ and any nonempty open set $V$ (not necessarily an open neighborhood of $x$) in $X$ there is a point $y$ in $V$ such that $\limsup_{n \to \infty} d(f^n(x), f^n(y)) \ge \delta$ and $\liminf_{n \to \infty} d(f^n(x), f^n(y)) = 0$.
\end{abstract}

\section{Introduction}
Let $\Sigma_2 = \{ \be : \be = \be_0\be_1 \cdots$, where $\be_i = 0$ or 1 $\}$ be the metric space with metric $d$ defined by $d(\be_0\be_1 \cdots$, $\ga_0\ga_1 \cdots) = \sum_{i=0}^\iy  {|\be_i-\ga_i|}/{2^{i+1}}$ and let $\si$ be the shift map defined by $\si(\be_0\be_1 \cdots) = \be_1\be_2 \cdots$.  The shift map $\si$ is often used {\bf{\cite{dev}}}, {\bf{\cite{rob}}} to model the chaoticity of a dynamical system.   But, what is the shift map chaotic about?     It is well known that the shift map has points with dense orbits (and hence is topologically transitive), has dense periodic points and has sensitive dependence on initial conditions.  Sensitive dependence on initial conditions, which is easily understood intuitively as nearby points, however close, will eventually separate a distance, is generally believed to be the central ingredient of chaos.  However, does it really reveal the true nature of chaos?  In {\bf{\cite{banks}}}, it is shown that sensitive dependence on initial conditions is a consequence of topological transitivity and dense periodic points and hence is a topological property.  On the other hand, if we let $W$ denote the dense subset of $\Sigma_2$ which consists of all elements with finitely many $1$'s in its expansion, then it is easy to see that the shift map is topologically transitive on $W$ and has sensitive dependence on initial conditions.  Yet, every point of $W$ is eventually fixed, i.e., for every $\be$ in $W$, there is a positive integer $n$ such that $\si^n(\be) = \bar 0 = 000 \cdots$.  So, $W$ is a system no one would like to call it chaotic.  These seem to suggest that sensitive dependence on initial conditions tells only part, but not the whole, of the chaos story.  But then what is the other part?  We know that the shift map has a property called extreme sensitive dependence on initial conditions which is stronger than sensitive dependence on initial conditions.  That is, there exists a positive number $\dl$ (for the shift map, we can choose $\dl = 1$) such that for any $\al$ in $\Sigma_2$ and any open neighborhood $V$ of $\al$ there is a point $\be$ in $V$ such that $\limsup_{n \to \infty} d(\si^n(\al), \si^n(\be)) \ge \dl$ and $\liminf_{n \to \infty} d(\sigma^n(\al), \sigma^n(\be)) = 0$.  We also know {\bf{\cite{du1}}} that the shift map has a dense uncountable invariant 1-scrambled set ($S$ is a $\dl$-scrambled set for some $\dl > 0$ {\bf{\cite{ly}}} if and only if (i) for any $x \ne y$ in $S$, $\limsup_{n \to \infty} d(\si^n(x), \si^n(y)) \ge \dl$ and $\liminf_{n \to \infty} d(\si^n(x), \si^n(y)) = 0$ and (ii) for any $x$ in $S$ and any periodic point $p$ of $\si$, $\limsup_{n \to \infty} d(\si^n(x), \si^n(p)) \ge \dl/2$).  However, are extreme sensitive dependence on initial conditions and the existence of a dense uncountable invariant 1-scrambled set all that the shift map is chaotic about?  This motivates us to investigate the chaoticity of the shift map even further.  Surprisingly, we find that the shift map, although may not be more chaotic than we can imagine, is definitely more chaotic than we previously thought.

\section{The chaoticity of the shift map}
First we introduce some terminology.  For any two finite strings $D = \be_0\be_1 \cdots \be_m$ and $E = \ga_0\ga_1 \cdots \ga_n$ of $0$'s and $1$'s, let $B \cdot E = \be_0\be_1 \cdots \be_m \ga_0\ga_1 \cdots \ga_n$ denote the  concatenation of $D$ and $E$.  Sometimes, we simply write $DE$ for $D \cdot E$ when no confusion arises.  If $D_1, D_2, \cdots, D_k$ are $k$ finite strings of $0$'s and $1$'s, the concatenation $D_1 \cdot D_2 \cdot \ldots \cdot D_k$ of $D_1, D_2, \cdots, D_k$ are defined similarly.  For any $\ga_k = 0$ or $1$, let $\ga'_k = 0$ if $\ga_k = 1$ and $\ga'_k = 1$ if $\ga_k = 0$.  For any finite string $E = \ga_k\ga_{k+1}\ga_{k+2} \cdots \ga_{k+n}$ of $0$'s and $1$'s, let $E'$ denote the finite string $\ga'_k\ga'_{k+1}\ga'_{k+2} \cdots \ga'_{k+n}$.  For any $\ga = \ga_0\ga_1\ga_2 \cdots$ in $\Si_2$ and any positive integers $i$ and $j$, let $C(\ga,i:j) = \ga_i\ga_{i+1} \cdots \ga_j$ and let $B(\ga,i,j) = C(\ga,i:i+(j!/2)-1) \,\, C(\ga,i+(j!/2)-1): i+(2\cdot j!/2-2)) \,\, C(\ga, i+(2 \cdot j!/2-2): i+(3 \cdot j!/2-3)) \, \cdots \linebreak$ 
$C(\ga, i+(j\cdot j!/2-j): i+((j+1)j!/2-j-1))$ which consists of $j+1$ blocks of finite portions of $\ga$ of length $j!/2$ each and the last element of each block coincides with the first element of the next block.  So, the length of $B(\ga,i,j)$ is $(j+1)!/2$.  On the other hand, let $\hat B(\ga,i,j) = B(\ga,i,j) \cdot B'(\ga,i+(j+1)!/2,j)$ be the concatenation of $B(\ga,i,j)$ and  $B'(\ga,i+(j+1)!/2,j)$.  So, the length of $\hat B(\ga,i,j)$ is $(j+1)!$.  

It is evident that any element of $\Si_2$ which contains every finite sequence of $0$'s and $1$'s is a transitive point of $\si$, i.e., point with dense orbit, and there are uncountably many of such points.  Let $\al = \al_0\al_1 \cdots$ be a fixed transitive point in $\Si_2$ and, for any integer $m \ge 5$ and any $\gamma = \gamma_0\gamma_1\gamma_2 \cdots$ in $\Si_2$, let $A(\gamma, m) = \al_0\al_1 \cdots \al_{m!-1}(\gamma_0)^{(m-1)!} (\gamma_1)^{(m-1)!} \cdots (\gamma_{m-1})^{(m-1)!}(01)^{(m-2)!}(0011)^{(m-2)!} \cdots (0^{m-1}1^{m-1})^{(m-2)!}$, where $0^1 = 0, 0^2 = 00, 0^3 = 000, (01)^2 = 0101, (0011)^3 = 0011 \, 0011 \, 0011$, and so on.  Note that $A(\gamma, m)$ has length $3 \cdot m!$.  Let $\{ \, x_i : x_i = x_{i,0}x_{i,1}x_{i,2} \cdots, i \ge 1 \}$ be any given countably infinite set in $\Si_2$.  Since in our construction of the scrambled set $S$, the first few terms of every element of $S$ are not important, for simplicity, we let $k \ge 5$ be any {\it fixed} integer and let $b_j$, $0 \le j \le k!-1$ be any {\it fixed} $k!$ numbers of $0$'s and $1$'s.  For any $\gamma = \gamma_0\gamma_1\gamma_2 \cdots$ in $\Si_2$, let 

$\tau_\gamma = b_0b_1\cdots b_{k!-1} A(\gamma,k) \cdot \qquad \qquad \qquad \qquad \qquad \qquad \qquad \qquad \qquad \qquad \qquad \qquad \qquad \qquad \qquad \qquad \qquad \linebreak
\hat B(x_1,4\cdot k!, k-1) \hat B(x_2,4\cdot k!+k!, k-1) \hat B(x_3,4\cdot k!+2 \cdot k!, k-1) \cdots \hat B(x_{k-3},4\cdot k!+(k-4) \cdot k!, k-1) \cdot \linebreak
A(\gamma,k+1) \cdot \hat B(x_1,4\cdot (k+1)!, k) \hat B(x_2,4\cdot (k+1)!+(k+1)!, k) \cdots \hat B(x_{k-2},4\cdot (k+1)!+(k-3) \cdot (k+1)!, k) \cdot \linebreak$
$A(\gamma,k+2) \cdot \hat B(x_1,4(k+2)!,k+1) \hat B(x_2,4(k+2)!+(k+2)!,k+1) \cdots \hat B(x_{k-1},4(k+2)!+(k-2)(k+2)!,k+1)\cdots$.

Note that the string $b_0b_1\cdots b_{k!-1}$ has length $k!$, the string $b_0b_1\cdots b_{k!-1} A(\gamma,k)\cdot \hat B(x_1,4\cdot k!, k-1)$ $\cdots \hat B(x_{k-3},4\cdot k!+(k-4) \cdot k!, k-1)$ has length $(k+1)!$, and so on.  Therefore, for every $m \ge k$, $\si^{m!}(\tau_\gamma) = A(\gamma, m) \cdot \hat B(x_1, 4\cdot m!, m-1) \cdots$.  Let $S = \{ \sigma^n(\tau_\gamma): n \ge 0, \gamma \in \Si_2 \}$.  Then it is clear that $S$ is a dense uncountable invariant (i.e., $\si(S) \subset S$) subset of $\Si_2$ which consists of transitive points.  We now use the blocks $A(\gamma, m)$'s in the definition of $\tau_\gamma$ to show that $S$ is a 1-scrambled set for $\sigma$.  Let $\ga = \ga_0\ga_1 \cdots$ and $\be = \be_0\be_2 \cdots$ be distinct points in $\Si_2$ with $\ga_s \ne \be_s$ and let $0 \le i < j$.  In the definition of $\tau_\gamma$, we examine closely the block $A(\gamma, m) =$ 
\qquad $\al_0\al_1 \cdots \al_{m!-1}(\gamma_0)^{(m-1)!} (\gamma_1)^{(m-1)!} \cdots (\gamma_{m-1})^{(m-1)!}(01)^{(m-2)!}(0011)^{(m-2)!} \cdots (0^{m-1}1^{m-1})^{(m-2)!}.$
  
In the block $A(\gamma, m)$ with $m \ge k$, because we have the string $(\gamma_s)^{(m-1)!}(\gamma_{s+1})^{(m-1)!}\cdots$ of $0$'s and $1$'s, we obtain that $\si^{2m!+s \cdot (m-1)!}(\tau_\gamma) = (\gamma_s)^{(m-1)!}(\gamma_{s+1})^{(m-1)!}\cdots$.  Similarly, $\si^{2m!+s \cdot (m-1)!}(\tau_\beta) = (\beta_s)^{(m-1)!}$ $(\beta_{s+1})^{(m-1)!}\cdots$.  So, $\limsup_{n \to \iy} d(\si^n(\si^i(\tau_\beta), \si^n(\si^i(\tau_\gamma)) \ge \lim_{m \to \iy}$ $d((\beta_s)^{(m-1)!}\cdots, (\gamma_s)^{(m-1)!}\cdots) \ge \lim_{m \to \iy} 1/2 + 1/2^2 + 1/2^3 + \cdots + 1/2^{(m-1)!} = 1$.  On the other hand, let $t_m = 3m!+(j-i)(j-i-1)\cdot (m-2)!$.  Then, because $A(\gamma, m)$ has the finite string $(0^{j-i}1^{j-i})^{(m-2)!}$ of $0$'s and $1$'s, although $\gamma \ne \beta$, we still have $\si^{t_m}(\tau_\gamma) = (0^{j-i}1^{j-i})^{(m-2)!} \cdots = 0^{j-i}(1^{j-i}0^{j-i})^{(m-2)!-1}1^{j-i} \cdots$ and $\si^{t_m}(\tau_\beta) = 0^{j-i}(1^{j-i}0^{j-i})^{(m-2)!-1}1^{j-i} \cdots$.  Thus, $\si^{t_m+j-i}(\tau_\gamma) = (1^{j-i}0^{j-i})^{(m-2)!-1}1^{j-i} \cdots = 1^{j-i}(0^{j-i}1^{j-i})^{(m-2)!-1} \cdots$.  Consequently, $\limsup_{n \to \iy}$ $d(\si^n(\si^i(\tau_\beta))$, $\si^n(\si^j(\tau_\gamma))) \ge$ $\lim_{m \to \iy} d(\si^{t_m}(\tau_\beta)$, $\si^{t_m+j-i}(\tau_\gamma))$ $\ge \lim_{m \to \iy}$ $1/2 + 1/2^2 + \cdots + 1/2^{j-i+2\cdot(j-i)(m-2)!} + \cdots = 1$.  Therefore, for any $x \ne y$ in $S$, we have $\limsup_{n \to \iy}$ $d(\si^n(x), \si^n(y))$ $= 1$.  Finally, since in the definition of $\tau_\beta$ and $\tau_\gamma$, there are infinitely many $n_i$'s such that both $\si^{n_i}(\tau_\beta)$ and $\si^{n_i}(\tau_\ga)$ start with the same arbitrarily long strings of 0's, we obtain that $\liminf_{n \ra \iy} d(\si^n(x), \si^n(y)) = 0$ for any $x \ne y$ in $S$ and similarly, since $\tau_\gamma$ contains arbitrarily long string of 0's, we obtain that $\limsup_{n \ra \iy} d(\si^n(x)$, $\si^n(y)) \ge 1/2$ for any $x$ in $S$ and any periodic point $y$ of $\si$.  This shows that $S$ is a dense uncountable invariant 1-scrambled set for $\sigma$.

Now let $x_i = x_{i,0}x_{i,1}x_{i,2} \cdots$ and let $\si^j(\tau_\gamma)$ be a point in $S$.  For any $m > k + i + j + 3$, let $t(m,i,j) = 4m!+(i-1)\cdot m!+j \cdot (m-1)!/2$.  Then $\si^{t(m,i,j)}(\si^j(\tau_\gamma)) = x_{i,t(m,i,j)}x_{i,t(m,i,j)+1} \cdots x_{i,t(m,i,j)+(m-1)!/2-j} \cdots$ and $\si^{t(m,i,j)}(x_i) = x_{i,t(m,i,j)}x_{i,t(m,i,j)+1} \cdots x_{i,t(m,i,j)+(m-1)!/2-j} \cdots$.  Similarly, $\si^{t(m,i,j)+m!/2}(\si^j(\tau_\gamma)) = x_{i,t(m,i,j)+m!/2}'x_{i,t(m,i,j)+m!/2+1}' \cdots x_{i,t(m,i,j)+m!/2+(m-1)!/2-j}' \cdots$ and $\si^{t(m,i,j)+m!/2}(x_i) =$ $x_{i,t(m,i,j)+m!/2} \linebreak$ $x_{i,t(m,i,j)+m!/2+1} \cdots x_{i,t(m,i,j)+m!/2+(m-1)!/2-j} \cdots$.  Therefore, we easily obtain that $\limsup_{n \to \iy} d(\si^n(x_i)$, $\si^n(\si^j(\tau_\gamma))) = 1$ and $\liminf_{n \to \iy} d(\si^n(x_i), \si^n(\si^j(\tau_\gamma))) = 0$.  This proves the following result.
  
{\bf Theorem 1.}
{\it For any given countably infinite subset $X$ of $\Sigma_2$, there exists a dense uncountable invariant 1-scrambled set $Y$ of transitive points in $\Sigma_2$ such that, for every $x$ in $X$ and every $y$ in $Y$, $\limsup_{n \ra \iy} d(\si^n(x), \si^n(y)) = 1$ and $\liminf_{n \ra \iy}$ $d(\si^n(x), \si^n(y)) = 0$.}

For the tent map $T(x) = 1 - |2x -1|$ on $[0, 1]$, we can use the symbolic representations of points in $[0, 1]$ as introduced in {\bf{\cite{du1}}} and use a similar variant of the above $\tau_\gamma$ to show that a result similar to Theorem 1 also holds for $T$.

{\bf Theorem 2.}
{\it For any given countably infinite subset $X$ of $[0, 1]$, there exists a dense uncountable invariant 1-scrambled set $Y$ of transitive points in $[0, 1]$ such that, for every $x$ in $X$ and every $y$ in $Y$, $\limsup_{n \ra \iy} |T^n(x) - T^n(y)| \ge 1/2$ and $\liminf_{n \ra \iy}$ $|T^n(x) - T^n(y)| = 0$.}

\section{The true nature of chaos}
Theorem 1 has a very important consequence.  That is, given any point $x$ in $\Sigma_2$ then at just about everywhere (the corresponding dense set $Y$) in $\Sigma_2$, whether it is close to  $x$ or far away from it we can always find a point $y$ (in $Y$) whose iterates satisfy $\limsup_{n \ra \iy} d(\si^n(x), \si^n(y)) = 1$ and $\liminf_{n \ra \iy}$ $d(\si^n(x), \si^n(y)) = 0$.  This seems to suggest that in a chaotic system not only nearby points can separate apart but also far apart points can get close to each other and these happen infinitely often.  After all, who knows when is the very beginning time in this ever changing world?  Two present far apart points may be very close to each other some time earlier.  

Theorem 1 also reveals a very striking property for the shift map.  That is, when we let $X = \{ x_0, \sigma(x_0), \sigma^2(x_0), \cdots \}$ denote the orbit of any given point $x_0$ in $\Sigma_2$, then Theorem 1 implies the existence of a dense uncountable invariant 1-scrambled set $Y$ of transitive points in $\Sigma_2$ such that, for every positive integer $m$ and every $y$ in $Y$, $\limsup_{n \ra \iy} d(\si^n(\si^m(x_0)), \si^n(y)) = 1$ and $\liminf_{n \ra \iy}$ $d(\si^n(\si^m(x_0)), \si^n(y)) = 0$.  In particular, this says that, for any point $x_0$ and any time $m$ earlier, at about everywhere (the corresponding dense set $Y$) in $\Sigma_2$, we can find a point $y$ (in $Y$) whose trajectory eventually catches up with that of the point $x_0$ to within any prescribed distance (since $\liminf_{n \ra \iy}$ $d(\si^n(\si^m(x_0)), \si^n(y)) = 0$) even though $x_0$ starts out time $m$ earlier than $y$.  

\section{A definition of chaos}
Let $(X, d)$ be an infinite compact metric space with metric $d$ and let $f$ be a continuous map from $X$ into itself.  We say that $f$ is chaotic (cf. {\bf{\cite{ak}}}, {\bf{\cite{ko}}}) if there exists a positive number $\delta$ such that for any point $x$ and any nonempty open set $V$ (not necessarily an open neighborhood of $x$) in $X$ there is a point $y$ in $V$ such that $\limsup_{n \ra \iy} d(f^n(x), f^n(y)) \ge \delta$ and $\liminf_{n \ra \iy} d(f^n(x), f^n(y)) = 0$.  Our definition of chaos is stronger than that of Li-Yorke sensitivity considered in {\bf{\cite{ak}}}.  In Theorem 4, we give an example which is Li-Yorke sensitive but not chaotic. By Theorems 1 and 2 and {\bf{\cite{au}}}, {\bf{\cite{ze}}}, we have the following result.

{\bf Theorem 3.}
{\it The following statements hold.
\begin{itemize}
\item[(a)] 
The shift map $\sigma$ is chaotic on $\Sigma_2$.

\item[(b)]
The tent map $T(x) = 1 - |2x -1|$ is chaotic on $[0, 1]$.

\item[(c)]
$F_\mu$ is chaotic on $\Ld_\mu$ for any $\mu \ge 4$, where $F_\mu(x) = \mu x(1-x)$ and $\Ld_\mu = \cap_{n=0}^\iy F_\mu^{-n}([0, 1])$ for $\mu > 4$ and $F_\mu(x) = [0, 1]$ for $\mu = 4$.
\end{itemize}}

The chaotic maps in Theorem 3 are all topologically transitive and Li-Yorke sensitive.  However, not all topologically transitive maps or Li-Yorke sensitive maps are chaotic.  The following is such an example.  

{\bf Theorem 4.}
{\it Let $g(x)$ be a continuous map from $[-1,1]$ onto itself defined by letting $g(x) = 2x+2$ for $-1 \le x \le -1/2$; $g(x) = -2x$ for $-1/2 \le x \le 0$; and $g(x) = -x$ for $0 \le x \le 1$.  Then $g$ is topologically transitive and Li-Yorke sensitive but not chaotic because the period-2 point $-2/3$ and the interval $[0, 1]$ are juming alternatively and never get close to each other.}

Chaotic maps are not necessarily topologically transitive.  The following is such an example.  

{\bf Theorem 5.}
{\it Let $T(x) = 1 - |2x - 1|$ for $0 \le x \le 1$ and let $h$ be a continuous map from $[-1/2, 1]$ to itself defined by $h(x) = -x$ for $-1/2 \le x \le 0$ and $h(x) = T(x)$ for $0 \le x \le 1$.  Then $h$ is chaotic on $[-1/2, 1]$ but $h$ is not topologically transitive.}

If $f$ is chaotic on the compact interval $I$ in the real line, then $f^2$ is turbulent (i.e., there exist closed subintervals $I_0$ and $I_1$ of $I$ with at most one point in common such that $f^2(I_0) \cap f^2(I_1) \supp I_0 \cup I_1$) as shown below.

{\bf Theorem 6.}
{\it Let $I$ denote a compact interval in the real line and let $f$ be a continuous map from $I$ into itself.  If $f$ is chaotic, then $f^2$ is turbulent.}

\noindent
{\it Proof.}
Let $z$ be a fixed point of $f$ and let $W$ be a nonempty open set in $I$.  Assume that $f$ is chaotic.  Then there exist a positive number $\delta$ and a point $c$ in $W$ such that $\limsup_{n \ra \iy} |f^n(c) - f^n(z)| \ge \delta$ and $\liminf_{n \ra \iy}$ $|f^n(c) - f^n(z)| = 0$.  In particular, the $\omega$-limit set $\omega(f, c)$ of $c$ with respect to $f$ ($x$ is in $\omega(f ,c)$ if and only if $\lim_{k \to \infty} f^{n_k}(c) = x$ for some sequence of positive integers $n_k \to \infty$) contains the fixed point $z$ of $f$ and a point different from $z$.  Therefore, $f^2$ is turbulent.
\hfill\sq

\end{document}